\documentclass{amsproc}
\usepackage{amsmath,amsfonts,amssymb,amsthm,epsfig,color,tikz}

\newcommand{\bb}{\mathbb}

\newcommand{\x}{\mathbf{x}}
\newcommand{\C}{\bb C}

\newcommand{\Z}{\bb Z}
\newcommand{\R}{\bb R}
\newcommand{\N}{\bb N}

\newcommand{\om}{\omega}

\newtheorem*{lemma*}{Lemma}
\newtheorem*{question*}{Question}

\newtheorem*{theorem*}{Theorem}

\numberwithin{Def}{section}

\newcommand{\vv}{\mathbf{v}}

\newcommand{\hh}{\mathcal H}

\newcommand{\La}{\Lambda}
\newtheorem{Theorem}{Theorem}

\newtheorem{lemma}[Theorem]{Lemma}

\begin{document}
\title{Random Affine Lattices}
\author{Jayadev~S.~Athreya}
\email{jathreya@iilinois.edu}
\address{Department of Mathematics, University of Illinois Urbana-Champaign, Urbana, IL}
\thanks{J.S.A partially supported by NSF grant
   DMS 1069153 and grants DMS 1107452, 1107263, 1107367  ``RNMS: GEometric structures And Representation varieties" (the GEAR Network)."} 
   
   \begin{abstract} We give a bound on the probability that a randomly chosen affine unimodular lattice has large holes, and a similar bound on the probability of large holes in the spectrum of a random flat torus. We discuss various motivations and generalizations, and state several open questions.
   \end{abstract}
\maketitle

\begin{center} \it In honor of Professor Ravi Kulkarni \end{center}

\subsection{Introduction and motivation}\label{sec:intro} In this note, we study the properties of \emph{random affine unimodular lattices} and in Euclidean spaces, and associated problems on the spectra of random flat metrics on tori. The properties of random geometric structures, random surfaces, and random flat tori have been of great interest to mathematicians (and physicists) in many different fields, including number theory, geometry, and probability theory. 

Many different models of random surfaces have been proposed for different purposes: for example, Sheffield's extraordinary thesis~\cite{Sheffield} describes some of the beautiful work on random surfaces arising from random functions on the integer lattice studied by, among others, Cohn, Elkies, and Propp~\cite{CEP} and Kenyon, Okounkov, and Sheffield~\cite{KOS}. For Riemann surfaces and hyperbolic surfaces, models have been developed and studied by, for example, Brooks-Makover~\cite{BM} and Mirzakhani~\cite{Mirzakhani}.

Motivations for studying random flat tori come from wanting to understand various geometric and analytic properties, in particular, counting problems and associated problems on the spectrum of the Laplacian. Problems such as the Berry-Tabor conjecture and related problems on eigenvalue spacings of (random) flat tori have been studied by, among many others, Eskin-Margulis-Mozes~\cite{EMM}, Petridis~\cite{Petridis} and VanDerKam~\cite{VanDerKam} (see~\cite{Marklof} for an excellent survey on the Berry-Tabor conjecture). 

Associated to each unit-volume $d$-dimensional flat torus $T$ is the unimodular lattice $\La \subset \R^d$ so that $T = \R^d/\La$. Motivated by problems on dynamics of unipotent flows, the author and G.~Margulis~\cite{AMarg} studied the properties of random unimodular lattices, and their estimates were extended and sharpened in certain cases by Strombergsson~\cite{Strom}, who was motivated his joint work with Marklof (see, for example,~\cite{MS}) on dynamics in the periodic Lorentz gas. The main results in our paper are to extend the work on random unimodular lattices in~\cite{AMarg} to the setting of random affine lattices, and to give applications of~\cite{AMarg} to the structure of the Laplacian for a random flat torus.

\subsection{The space of lattices}\label{sec:lattices}  An \emph{unimodular lattice} in $\R^d$ is a discrete additive subgroup of covolume $1$, and an \emph{affine unimodular lattice} in $\R^d$ is the orbit of a point $\x \in \R^d$ under translations by a unimodular lattice.  Thus, any affine unimodular lattice $\La$ has the form $$\La = g\Z^d + \x,$$ where $g \in SL(d, \R)$ and $\x$ is well-defined up to translation by $g \Z^d$, that is, it can be viewed as an element of the $d$-torus $\R^d/g\Z^d$. We can parameterize the space $Y_d$ of unimodular affine lattices by the homogeneous space $SL(d, \R) \ltimes \R^d/SL(d, \Z) \ltimes \Z^d$, via the identification $$\La = g\Z^d + \x \longmapsto (g, \x) SL(d, \Z) \ltimes \Z^d.$$ We can endow $Y_d$ with the probability measure $\mu = \mu_d$ induced by Haar measure on $SL(d,\R) \ltimes \R^d$. Our first main result, Theorem~\ref{theorem:affmink}, is on the probability that a lattice chosen from $Y_d$ at random according to $\mu_d$ has large `holes'.

\begin{Theorem}\label{theorem:affmink} Let $A \subset \R^d$ be a measurable set. Then $$\mu(\La \in Y_d: \La \cap A = \emptyset ) < \frac{1}{1+|A|},$$ where $|A|$ denotes the Lebesgue measure of $A$.
\end{Theorem}

\noindent For the space $X_d = SL(d, \R)/SL(d, \Z)$ of unimodular lattices,  the analogous result (with respect to $\nu = \nu_d$, the Haar probability measure on $X_d$) was shown by the author and G.~A.~Margulis~\cite[Theorem 2.2]{AMarg}, and was described as a `Random Minkowski Theorem', to emphasize the connection to the classical Minkowski convex body theorem. The bound in the setting of regular lattices was of the form $\frac{C_d}{|A|}$, where the constant $C_d$ could be taken to be $8 \zeta(d-1)/\zeta(d)$, for $d \geq 3$, and $16\zeta(2)$ for $d=2$. Here, the constant $1$ is independent of dimension. The order of approximation, inversely proportional to the volume of $|A|$ was shown to be sharp in the setting of regular lattices by Strombergsson~\cite{Strom}, even if $A$ was restricted to be convex. It would be interesting to understand the sharp values of these constants for both the affine and regular lattice settings.

Our second main result, Theorem~\ref{theorem:spectra} is on the spectrum of a random flat $d$-dimensional torus. Given an unimodular lattice $\La \subset \R^d$, let $\Omega(\La) \subset \R^+$ be the spectrum of the Laplacian on the flat torus $\R^d/\La$. This consists of lengths of (non-zero) vectors in the dual lattice $\La^*$. Let $C_d$ be as in the above discussion, and let $N_d$ denote the volume of the unit ball in $\R^d$. Then we have

\begin{Theorem}\label{theorem:spectra} Let $A \subset \R^+$ be a measurable set, and let $|A|_d = d N_d\int_{A} r^{d-1} dr$. Then $$\nu (\La \in X_d: \Omega(\La) \cap A = \emptyset) \le \frac{C_d}{|A|_d}.$$
\end{Theorem}

As discussed in \S\ref{sec:intro}, the structure of the spectrum $\Omega(\La)$ has been extensively studied. For example, Petridis~\cite{Petridis}, showed that the (appropriately normalized) spacing set of a random flat torus is dense in $\R^d$ with probability $1$ (a result that neither implies nor is implied by our result). VanDerKam~\cite{VanDerKam} showed that with probability $1$, the eigenvalues of a random flat torus in dimension $d$ follow a uniform distribution with respect to their $1, \ldots, \lfloor \frac d 2 \rfloor$-spacings. Eskin-Margulis-Mozes~\cite{EMM} showed that in the setting of two-tori, the set of exceptions to the Berry-Tabor conjecture on consecutive eigenvalue spacings is in fact of Hausdorff dimension $0$ (more than just measure $0$). 

\subsection{Mean and Variance}\label{sec:mean} The proof of Theorem~\ref{theorem:affmink} is quite similar to the proof of the Random Minkowski Theorem given in~\cite{AMarg} for regular lattices, but is in some ways more transparent, especially in the setting $d=2$. Our key tools will be lemmas calculating the mean and variance of the random variable $\widehat{\chi}_A: Y_d \rightarrow \N \cup \{0\}$, given by $$\widehat{\chi}_A (\La) = \#( \La \cap A),$$ which counts the number of lattice points in $A$. These estimates allow us to use Chebyshev's theorem to estimate the probability that $\widehat{\chi}_A =0$. 

 We state our key lemmas and prove Theorem~\ref{theorem:affmink}. Fix $A \subset \R^d$ measurable, and as above define $\widehat{\chi}_A: Y_d \rightarrow \N \cup \{0\}$ as the counting function of the number of lattice points in $A$. We want to understand $$p_A : = \mu( \La \in Y_d: \widehat{\chi}_A(\La) =0).$$ By abuse of notation, let $$\mu_A : = \int_{Y_d} \widehat{\chi}_A d\mu_d$$ and $$\sigma^2_A:= \left(\int_{Y_d} \widehat{\chi}^2_A d\mu_d\right) - \mu^2_A$$Chebyshev's inequality yields that $$p_A \le \frac{\sigma^2_A}{\sigma^2_A + \mu^2_A}$$

\subsubsection{Mean} The following is the affine version of the classical \emph{Siegel integral formula} (see, for example, \cite[(3.6)]{ElBaz}):

\begin{lemma}\label{lem:mean} $\mu_A = |A|$.
\end{lemma}

\noindent The philosophy of the proof of this result is simple. Given $f \in C_c(\R^d)$, consider the functional $$f \longmapsto \int_{Y_d} \left(\sum_{v \in \La} f(v)\right) d\mu_d(\La).$$ Assuming this is well-defined, it is a linear $SL(d, \R) \ltimes \R^d$-invariant functional on $C_c(\R^d)$, and thus must come from integration against an $SL(d, \R) \ltimes \R^d$-invariant measure on $\R^d$, namely, a multiple of Lebesgue measure. By considering the indicator function of a large ball and standard lattice point counting arguments, the normalizing constant must be $1$, and thus, we have $$ \int_{Y_d}\left( \sum_{v \in \La} f(v) \right)d\mu_d(\La) = \int_{\R^d} f.$$ Applying this to $f = \chi_A$ we obtain the lemma. If $A$ is not precompact, a standard approximation argument gives the result. The assumption that the functional (here and below) is well-defined is justified by the well-definedness in the classical (i.e. non-affine) setting, since the space of affine unimodular lattices is a compact fiber bundle over the space of unimodular lattices.

\subsubsection{Variance} Similarly, we can modify a theorem of Rogers~\cite{Rogers} (see, for, example, \cite[(3.7)]{ElBaz}) to obtain:

\begin{lemma}\label{lem:var} $\sigma^2_A = |A|$
\end{lemma}

\noindent The philosophy of proof is similar. Given $h \in C_c(\R^d \times \R^d)$, consider the functional $$h \longmapsto \int_{Y_d} \left(\sum_{v_1, v_2 \in \La} h(v_1, v_2)\right) d\mu_d(\La)= \int_{Y_d} \left(\sum_{v_1, v_2 \in \La, v_1 \neq v_2} h(v_1, v_2)\right) d\mu_d(\La)  + \int_{Y_d} \left(\sum_{v \in \La} h(v, v)\right) d\mu_d (\La).$$ As above, this is well-defined, so it is a linear $SL(d, \R) \ltimes \R^d$-invariant functional on $C_c(\R^d \times \R^d)$ (where the action is the diagonal action), and thus comes from an $SL(d, \R) \times \R^d$-invariant measure on $\R^d \times \R^d$. The orbits of this action (for $d \geq 2$) consist of pairs of distinct vectors and pairs of equal vectors (corresponding to the decomposition of the sum above). Taking careful care of normalization, this yields $$\int_{Y_d} \left(\sum_{v_1, v_2 \in \La} h(v_1, v_2)\right) d\mu_d(\La) = \int_{\R^d \times \R^d} h(v_1, v_2) + \int_{\R^d} h(v, v).$$ Applying this to $h(v_1, v_2) = \chi_A(v_1) \chi_A(v_2)$, we obtain $$\int_{Y_d} \widehat{\chi}^2_A d\mu_d = |A|^2 + |A|.$$ Since $\mu_A = |A|$, we have $$\sigma^2_A = |A|^2 + |A| - |A|^2 = |A|$$ 

\subsubsection{Regular lattices}\label{sec:reg} In the $d=2$ setting, this proof is considerably simpler than the corresponding proof for regular (that is, non-affine) unimodular lattices. For regular unimodular lattices, this philosophy of proof was developed by Rogers~\cite{Rogers}, and requires us to understand the orbits of $SL(d, \R)$ on $\R^d \times \R^d$. For $d \geq 3$, a similar analysis to above yields the classical Random Minkowski theorem~\cite[Theorem 2.2]{AMarg}. However, for $d=2$, the wedge product (or determinant) is an invariant, and so the orbit structure is quite complicated, and this case required a very different analysis.

\subsubsection{Proof of Theorem~\ref{theorem:affmink}}\label{sec:proof}
By Chebyshev's inequality we have $$p_A \le \frac{\sigma^2_A}{\sigma^2_A + \mu^2_A} = \frac{|A|}{|A|+|A|^2} = \frac{1}{1+|A|}$$\qed\medskip

\subsection{Spectra} To prove Theorem~\ref{theorem:spectra}, note that $$\nu (\La \in X_d: \Omega(\La) \cap A = \emptyset) = \nu(\La^* \in X_d: \La^* \cap \tilde{A} = \emptyset),$$ where $\tilde{A} = \left\{ v \in \R^d: \|v\| \in A\right\}$. By~\cite[Theorem 2.2]{AMarg},  $$\nu(\La^* \in X_d: \La^* \cap \tilde{A} = \emptyset) \le \frac{C_d}{|\tilde{A}|}.$$ We have $$|\tilde{A}| = d N_d\int_{A} r^{d-1} dr,$$ yielding the result.\qed\medskip

\subsection{Open problems} In this section, we discuss several related open problems.

\subsubsection{Sharp constants} As discussed in \S\ref{sec:lattices}, it would be very interesting to get optimal values of constants for the upper bounds in Theorem~\ref{theorem:affmink} and~\cite[Theorem 2.2.]{AMarg} (which would in turn give information about Theorem~\ref{theorem:spectra}). That is, we would like to compute the quantities $$\limsup_{|A| \rightarrow \infty} (1+|A|) \mu(\La \in Y_d: \La \cap A = \emptyset)$$ and $$\limsup_{|A| \rightarrow \infty} |A| \mu(\La \in X_d: \La \cap A = \emptyset).$$

\subsubsection{Elementary proof for $n=2$} As mentioned in \S\ref{sec:reg}, the current proof of~\cite[Theorem 2.2]{AMarg} in dimension $d=2$ requires a different method than the one discussed in this paper. It uses a subtle analysis involving Eisenstein series and the representation theory of $SL(2, \R)$. It would be of great interest to give a more elementary proof in this setting, which in some sense is the simplest of all settings.

\subsubsection{Torsion affine lattices} There are other interesting $SL(d, \R)$-invariant probability measures on $Y_d$, supported on the set of affine unimodular lattices whose translation vector is torsion in the associated flat torus. It would be nice to investigate the probabilities that random lattices chosen according to these measures have large holes.

\subsubsection{Random translation surfaces} A natural family of generalizations of flat $2$-tori are given by \emph{translation surfaces}: A translation surface is a pair $(M, \om)$, where $M$ is a Riemann surface and $\om$ a holomorphic $1$-form. Each flat $2$-torus (together with the one-form $dz$) is naturally a translation surface. By fixing the combinatorics of the differential, we obtain spaces of translation surfaces known as \emph{strata}. Each stratum $\hh$ carries a natural absolutely continuous probability measure, known as \emph{Masur-Veech} measure $\mu_{MV}$. We refer the reader to~\cite{Zorich} for an excellent survey on translation surfaces

A \emph{saddle connection} on a translation surface $(M, \omega)$ is a geodesic (in the flat metric induced by $\om$) $\gamma$ connecting two zeros of $\omega$ (with none in its interior). To each saddle connection $\gamma$ one can associate a holonomy vector $\vv_{\gamma} = \int_{\gamma} \om \in \C$, by integrating the one-form. The set of holonomy vectors $\La_{sc}(\om)$ is a discrete subset of $\C \cong \R^2$. For flat $2$-tori, if we mark the point $0$ as a `singular' point, the set of saddle connections exactly corresponds to the set of \emph{primitive} vectors in the associated lattice. A natural question then becomes: can we find a constant $c = c_{\hh}$ associated to each stratum $\hh$ so that $$\mu_{MV}(\om \in \hh: \La_{sc}(\om) \cap A) \le \frac{c}{|A|}$$ for any $A \subset \R^2$?

\end{document}